\documentclass[12pt]{article}
\usepackage{amssymb}
\usepackage{latexsym}

\title{Harmonic tori and generalised Jacobi varieties.}
\author{Ian McIntosh\thanks{Partially supported by NSF grants DMS-9626804 
and DMS-9705479}\\
Dept.\ of Mathematics\\
University of York\\
Heslington, York, YO10 5DD, UK}
\date{December, 1999}

\newcommand{\C}{\mathbb{C}}
\newcommand{\R}{\mathbb{R}}
\renewcommand{\P}{\mathbb{P}}
\newcommand{\CP}{\mathbb{CP}}
\newcommand{\Z}{\mathbb{Z}}
\newcommand{\Tr}{\textrm{Tr}}

\newcommand{\Gr}{Gr_{k,n+1}}

\newcommand{\PU}{\mathbb{P}U_{n+1}}
\newcommand{\caL}{\mathcal{L}}
\newcommand{\caO}{\mathcal{O}}
\newcommand{\caE}{\mathcal{E}}

\newcommand{\caA}{\mathcal{A}}
\newcommand{\caU}{\mathcal{U}}
\newcommand{\caT}{\mathcal{T}}
\newcommand{\caH}{\mathcal{H}}
\newcommand{\fk}{\mathfrak{k}}
\newcommand{\fs}{\mathfrak{s}}
\newcommand{\fg}{\mathfrak{g}}
\newcommand{\fo}{\mathfrak{o}}
\newcommand{\fm}{\mathfrak{m}}

\newcommand{\Xo}{X^\prime}
\newtheorem{thm}{Theorem}
\newtheorem{prop}{Proposition}
\newtheorem{lem}{Lemma}
\newtheorem{cor}{Corollary}

\begin{document}
\maketitle

\section{Introduction.}

Over the last decade there has been considerable success in understanding
the construction of certain harmonic 2-tori in symmetric spaces (in particular,
the non-superminimal tori in $\C\P^n$ and $S^n$) using the methods of
integrable systems theory. For example, if $\varphi:M\to
S^2$ is a non-conformal harmonic torus one knows that $\varphi$ is determined
by its (real hyperelliptic) spectral curve $X$ equipped
with a degree two function $\lambda$. There are already 
several ways of describing the relationship between
the map $\varphi$ and its spectral data $(X,\lambda)$: see, for example,
\cite{Bob,Hit,PinS,McI1}. But we lack a direct, geometric picture
of how $\varphi$ arises from the algebraic geometry of $X$. 
What I want to present here is such a picture,
based on a remarkable property of certain generalised Jacobi varieties. Moreover, it 
is quite straightforward to extend this picture to produce
harmonic (indeed, pluri-harmonic) maps into Grassmannians and 
$\P U_{n+1}$ (i.e. $U_{n+1}/\hbox{
centre}$) and all of these will be maps of `finite type' in the sense of
\cite{BurFPP}. From this picture one also sees that some of these harmonic tori
are purely algebraic and these present an interesting class for further study.

To summarise, keeping with the example of non-conformal
tori in $S^2$, let us suppose for simplicity that $X$
is non-singular with Jacobi variety $J$. One of the necessary properties of
the pair $(X,\lambda)$ is that $\lambda$ is unbranched over $\lambda=1$: let
$O_1,O_2$ be the two points on $X$ with this value. Take $X^\prime$ to be the
singularisation of $X$ obtained by identifying $O_1$ with $O_2$ to obtain an
ordinary double point and let $J^\prime$ denote its generalised Jacobi variety. This is 
the moduli space of degree zero line bundles over $X^\prime$: it  is a
non-trivial group extension of $J$ by the multiplicative complex numbers $\C^\times$. 
One may think of a line bundle over $X^\prime$ as a 
line bundle over $X$ with an
identification between the fibres over $O_1$ and $O_2$ and it is this identification 
which is encoded by the fibres of $J^\prime\to J$.
I claim that the harmonic map $\varphi$ factors through
$J^\prime$ i.e.\ $\varphi=\psi\circ\gamma$ where 
\[
M\stackrel{\gamma}{\rightarrow}J^\prime\stackrel{\psi}{\rightarrow}\P^1.
\]
Moreover, the map $\psi$ is a rational map: it is the ratio of two 
sections of the
pullback to $J^\prime$ of the $\theta$-line bundle over $J$. This pullback has 
an infinite dimensional space of sections but there is a natural 
spanning set indexed by $\Z$. The sections we want have index $0$ and $1$. The map
$\gamma$ is a simply a homomorphism of real groups: we think
of $M$ as $\R^2$ modulo a lattice, while $J^\prime$ has a real structure
inherited from the real structure on $X$. The  subgroup of real points $J^\prime_R\subset
J^\prime$ is
a compact real torus of dimension $g+1$ for $g=\hbox{genus}(X)$.
Indeed, if $z=x+iy$ denotes the standard complex coordinate 
on $\R^2$, then $\gamma$ is completely determined by the property that
$\partial\gamma/\partial z$ is tangent to an
Abel map image of $X^\prime$ in $J^\prime$. Although $\psi$ is only rational it is
well-defined on the real slice $J^\prime_R$ of
$J^\prime$ in which $\gamma$ takes its image, hence $\psi\circ\gamma$ is
globally defined. Moreover, as we translate
$\gamma(M)$ around $J^\prime_R$ we obtain new harmonic maps: this is a well-known property of
maps of finite type. In
fact we will see that a translation along the fibres of $J^\prime_R\to J_R$ merely changes the
the image of $\varphi$ by an isometry.                       

The mechanism which underlies this construction is a beautiful fact about this
generalised Jacobian. Given a line bundle $\caL$ over $X$ of degree $g+1$ for which
$\caL(-O_1-O_2)$ is non-special one has, over a Zariski open subset of
$J$, a rank two vector bundle $E$ whose
fibre over $L\in J$ is $\Gamma(\caL\otimes
L)$, the space of globally holomorphic sections of $\caL\otimes L$. 
In our case this open set contains the real subgroup $J_R\subset J$ and, locally,
it is sections of $\P E$ (the projective bundle) 
over $J_R$ which provide the harmonic map (see \cite{McI1}). However, 
when one tries to extend these sections globally one encounters 
non-trivial holonomy around the periods of $J$. It
turns out that $J^\prime$ can be realised as a subbundle of the bundle of
projective frames of $E$. The holonomy of these sections of $E$ takes values in the non-trival
$\C^\times$-bundle $\pi:J^\prime\to J$. But $\pi^*J^\prime$ has
a tautological section and this section is algebraic, hence the bundle of projective 
frames of $\pi^*E$ is algebraically trivialisable. Therefore the harmonic tori factor
naturally through $J^\prime$.

Of course, to prove all this we need to have a working knowledge of generalised Jacobi varieties.
Since I do not assume many differential geometers will be familiar with this the first two sections
explain as much as we need to prove that $J^\prime$ is a natural subbundle of 
$\P Fr(E)$, the bundle of projective frames of $E$
(indeed this is done for the whole Picard group $Pic(X^\prime)$). It follows that $\P
Fr(\pi^*E)\to J^\prime$ has a canonical trivialisation. This is necessary to turn sections of
$Gr_k(\pi^*E)$ (the bundle of $k$-planes) and $\P Fr(\pi^*E)$ into maps from $J^\prime$ into
$\Gr$ and $\P Gl_{n+1}$. I also describe the appropriate Hermitian structure on fibres of
$\pi^*E$: this is necessary to be able to discuss
real maps into the symmetric spaces $\Gr$ and $\PU$.

It is section 3 which shows that, for the right choice of
sections restricted to the right choice of real two dimensional subgroup $\gamma(\R^2)\subset
J^\prime$, one obtains harmonic maps.
The most demanding part of the paper is to understand what constitutes `the right choice'.
Section 3 simply presents these choices as a \textit{fait accompli}, but the motivation is as
follows.

For simplicity let us return to the example of maps
$\varphi:\R^2\to S^2\simeq \CP^1$. It is well-known (see \cite{Uhl}) that the harmonic maps
are characterized by possessing a horizontal holomorphic lift into a loop group: we can
represent this by the extended frame $F_\lambda$, which maps $\R^2$ into 
$C^\omega(\C^\times,GL_2)$ (recall that $\varphi$ is recovered by taking the complex line
generated by the first column of $F_\lambda(z)$ evaluated at $\lambda=1$). 
We can think of the columns of $F_\lambda$ as
being smooth maps of $\R^2$ into the space $\caH$ of holomorphic sections of the trivial bundle
$\P^1_\lambda\backslash\{0,\infty\}\times\C^2$. When $\varphi$ is doubly periodic it was shown
in \cite{McI1} that $\caH$ is a rank 1 module over a commutative matrix algebra of symmetries
of $F_\lambda^{-1}dF_\lambda$ (these are the polynomial
Killing fields introduced by \cite{BurFPP}). This algebra gives us the double cover
$\lambda:X\to\P^1_\lambda$ and allows 
us to identify $\caH$ with the space of holomorphic sections of a holomorphic 
line bundle $\caL$ over $X\backslash\lambda^{-1}(\{0,\infty\})$. The columns of $F_\lambda(z)$
correspond to sections of $\caL$. However, to extend these sections globally over $X$ we must
make $\caL$ move with $z$. The motion can be explicitly determined from a 1-cocycle on $X$
determined by the behaviour of $F_\lambda^{-1}\partial F_\lambda/\partial z$ 
about $\lambda=0,\infty$.  One shows that $\caL(z)$ moves 
linearly along $J$ and tangent to an Abel image of $X$ in $J$, in the manner
described earlier. This gives us the right choice of subgroup $\gamma(\R^2)\subset J^\prime$. 
The complex line generated by the first column of $F_\lambda(z)$ corresponds to a line
in $\Gamma(\caL(z))\simeq\C^2$ determined purely by a divisor on $X$ 
supported over $\lambda=\infty$. The space $\P\Gamma(\caL(z))$ is the fibre of $Gr_1(\pi^*E)$
over $\caL(z)$ and we see that the right choice of section of this bundle is determined by a
divisor on $X$ over $\lambda=\infty$. 

The same principle is known to apply to all the (non-isotropic) harmonic tori in 
$\CP^n$ \cite{McI1,McI2}.
In section 3 I have adapted this principle to the construction of harmonic maps into $\Gr$
and $\PU$. 
The first steps are to describe the divisor over
$\lambda=0,\infty$, which picks out the right section of $Gr_k(\pi^*E)$ (and later,
$\P U(\pi^*E)$), and to describe the right
direction to move along $J^\prime$. Then $F_\lambda$ is constructed and shown to frame a
(pluri)-harmonic map.  I fully expect that this construction
gives all harmonic tori of semisimple finite type, in much the same way that one argues for $\CP^n$
\cite{McI1}. 
However, there is as yet no good description of the necessary and sufficient conditions for a torus
in a symmetric space  of rank greater than one to be of semisimple finite type, although some 
results have been obtained for $\Gr$ in \cite{Uda}.

The fourth section introduces the pullback of the $\theta$-line bundle to $J^\prime$ and shows 
that sections of this bundle provide maps $J^\prime\to\Gr$ (or $\PU$) which coincide with our harmonic maps
on the appropriate subgroup $\gamma(\R^2)\subset J^\prime$. 
We now see that the harmonic map will be an algebraic object whenever this subgroup 
$\gamma(\R^2)$ is algebraic i.e.\
whenever it is the real slice of an algebraic subgroup of complex dimension two.
I call these algebraic harmonic maps, although one must not be mislead: the (real) 
harmonic surface in $\Gr$ will
not necessarily be real algebraic (since the map $J^\prime\to Gr$ won't in general be real)
but it will be contained in a complex algebraic surface.
It is easiest to discuss the case of maps into $\CP^n$ and here I show that if 
$\gamma(\R^2)$ is algebraic it is
necessarily a two torus. Examples do exist: all $S^1\times S^1$-equivariant harmonic
tori in $\CP^n$ are algebraic (and have the Riemann sphere as their spectral curve); the Gauss maps
of the Delaunay surfaces are algebraic. But it remains
an open and interesting question to find other examples.

\smallskip\noindent
\textbf{Acknowledgements.} This paper was written while the author was visiting the Center for
Geometry, Analysis, Numerics and Graphics, Dept.\ of Mathematics, University of Massachusetts,
Amherst.

\subsection{Spectral data.}
\label{sec-spectral}

Throughout this article we will work with spectral data $(X,\lambda,\caL)$
where $X$ is 
a compact Riemann surface equipped with a rational function $\lambda$ of
degree $n+1$ and a line bundle $\caL$ over $X$ of degree $g+n$ where $g$ is
the genus of $X$. 
The aim is that from this data one obtains harmonic maps into $\Gr$ or $\P U_{n+1}$: certain
extra properties of this data distinguish which target is obtained and whether or not the map is
conformal. 
But in all cases $X$ must be equipped with real involution 
$\rho$ for which $\lambda\circ\rho = \bar\lambda^{-1}$,
$\lambda$ must have no branch points on the unit
circle and $\rho$ must fix all points over the pre-image of the this. Further, we insist that
$\caL$ satisfy the reality condition 
$\overline{\rho_*\caL}\simeq\caL^{-1}\otimes\caO_X(R)$ where $R$ is the ramification
divisor of $\lambda$ 
(it is the specification of $R$ which determines the target for the harmonic
map). 
Under these conditions one knows (cf.\  \cite{McI1,McI2}) that $\dim\Gamma(\caL)=n+1$ and
$\caL$ has no non-zero global sections vanishing at all the $n+1$ distinct points 
$O_1,\ldots,O_{n+1}$ over $\lambda=1$. Moreover we can 
use the trace map $\Tr:\Gamma(\caO_X(R))\to\Gamma(\caO_{\P^1})$ 
to equip each space $\Gamma(\caL)$ with a Hermitian inner product. To fix 
this inner product we require, first, a global section $\sigma$ so that $\Gamma(\caL)$ can be
realised as a divisor line bundle $\caO_X(D)$, and second, an isomorphism $f$ of
$\caO_X(D+\rho_*D)$ with $\caO_X(R)$ for which $\overline{\rho_*f}=f$ (we think of $f$ as a rational
function with divisor $D+\rho_*D-R$ and it can be chosen to be positive over the
unit circle $\vert\lambda\vert =1$, see \cite{McI2}). 
Then the Hermitian inner product on $\Gamma(\caL)$ is given by
\begin{equation}
\label{eq:h}
\begin{array}{ccc}
h(\sigma_1,\sigma_2)& =& \Tr(f(\sigma_1\sigma^{-1})\overline{\rho_*(\sigma_2\sigma^{-1})}\\
      & = & \sum_{j=1}^{n+1}f(O_j)
(\sigma_1\sigma^{-1})(O_j)\overline{\rho_*(\sigma_2\sigma^{-1})(O_j)},\end{array}
\end{equation}
Notice that if we change $\sigma$ then we must change $f$ too and as result we see that this Hermitian
inner product is unique up to a real scaling (the scale of $f$).

We will use $J$ to denote the Jacobi variety of $X$ and 
$J^\prime$ to denote the generalised Jacobi variety for the 
singularisation $X^\prime$ of $X$ obtained by identifying all the points $O_j$ together to 
obtain a single node $O\in X^\prime$. Thus $J^\prime$ is the Jacobian for the modulus 
$\fo=O_1+\ldots + O_{n+1}$, in the
terminology of Rosenlicht (see \cite{Ser}). 
The next section describes the relevant properties of $J^\prime$.

\smallskip\noindent
{\bf Notation.} If $E$ is a holomorphic vector bundle over a variety $M$ and $P\in M$ then $E_P$
will mean the stalk of the corresponding sheaf of locally holomorphic sections about $P$ while
$E\vert P$ will mean the fibre of $E$ over $P$ (and similar notation for restriction to subvarieties
of $M$).

\section{The generalised Jacobi variety.}

The aim in this section is to realise  the generalised Jacobian $J^\prime$ of $X^\prime$
as a subbundle of the projective bundle of frames of a particular rank $n+1$ vector bundle over $J$. 

Recall (from e.g.\ \cite{Ser}) that $J^\prime$ is, by definition, the group of (isomorphism
classes of) holomorphic
line bundles of degree zero over $X^\prime$. Equally, it is the the group of equivalence classes
of degree zero divisors on $X^\prime-\{O\}$ modulo the divisor $\fo$. Two such divisors
$D_1,D_2$ are equivalent modulo $\fo$ if $D_1-D_2$ is the divisor of a rational function
$f\in\C(X)$ which takes the same non-zero value 
at each point $O_j$. We will denote the divisor class of $D$ mod $\fo$ 
by $[D]_\fo$ and the divisor of any rational function $f$ by $(f)$.

Since this is a refinement of the divisor classes on $X$ the normalisation $X\to X^\prime$ induces,
by pullback, the extension
\begin{equation}
\label{eq:J}
1\to H_\fo\to J^\prime\stackrel{\pi}{\to} J\to 1
\end{equation}
where $H_\fo = \{[D]_\fo: D\hbox{  is the divisor of a rational function on $X$}\}$.
For our purposes it is more useful to describe $H_\fo$ as follows. Let $U_\fo$ be the
multiplicative group of rational functions on $X$ which have no poles or
zeroes on $\fo$ and let $U_O$ denote the subgroup of those which are obtained by
pullback from $\Xo$ (all such functions have $f(O_j)=f(O_1)$ for each point $O_j$).
Then $U_\fo/U_O$ is isomorphic to $H_\fo$: the
identification gives to each $f\in U_\fo$ its divisor class mod $\fo$. This
allows a local description of the fibres of $J^\prime$ over $J$. For each $L\in J$
let $L_\fo$ be the space of locally regular sections of $L$ about the support of $\fo$ and let
$L^\times_\fo$ be the subset of those which do
not vanish at any $O_j$. Then $U_O$ acts by multiplication on $L^\times_\fo$ and 
\begin{equation}
J^\prime\vert L\simeq L^\times_\fo/U_O.
\end{equation}
This equally well describes the fibres of the Picard variety extension
$\pi:Pic(X^\prime)\to Pic(X)$ where these are the groups of line bundles of any degree. 
Fom another point of view,
the fibre over $L\in Pic(X)$ consists of all n-tuples $(g_2,\ldots,g_{n+1})$ of linear bijections
$g_j:L\vert O_j\to L\vert O_1$.

Now we want to
realise $Pic(X^\prime)$ as a subbundle of a natural frame bundle over $Pic(X)$.
Let $\caE\to X\times Pic(X)$ be a 
Poincare line bundle i.e.\ a line bundle
with the property $\caE\vert X\times\{L\}\simeq L$. Recall (from e.g.\ \cite{ArbCGH})
that $\caE$ is not unique: if $\caT\to Pic(X)$ is any line bundle and $\pi_2:X\times Pic(X)\to
Pic(X)$ is projection on the second factor then $\caE\otimes\pi_2^*\caT$ is also a Poincare line
bundle. The normalisation
$X\to X^\prime$ pushes $\caE$ down to a coherent sheaf whose restriction to
$\{O \}\times Pic(X)$ we will call $E$. Its fibre over any point
$L$ is
\begin{equation}
\label{eq:E}
E\vert L\simeq L_\fo/(\fm_O\cdot L_\fo)= \oplus_{j=1}^{n+1} L\vert O_j
\end{equation}
where $\fm_O$ is the ideal in $\caO_\fo$ of locally regular functions about $\fo$ with
divisor of zeroes at least $\fo$. Thus $E$ is a rank $n+1$ vector bundle. Let $\P Fr(E)$ be the
principal bundle of projective frames of $E$ and let $Gr_k(\pi^*E)$ be the Grassmannian bundle of
$k$-planes in the fibres of $\pi^*E$. Since the effect of tensoring $\caE$ by $\pi_2^*\caT$ is to
globally rescale each line bundle $L$, these projective bundles are independent of the choice of
Poincare line bundle.
\begin{prop}
There exists a canonical injective morphism of principal bundles
$Pic(X^\prime)\to \P Fr(E)$. Consequently $\P Fr(\pi^*E)$ and $Gr_k(\pi^*E)$ are canonically
trivial.
\end{prop}
{\bf Proof.} Set $G^\C=\P GL_{n+1}$ and let $H\subset G^\C$ 
be the subgroup of diagonal matrices modulo
scaling. Consider $\P Fr(E)$ as a principal right $G^\C$-bundle over $J$ in the usual way.
In particular, for $f=[f_1,\ldots,f_n]$ a projective frame and $h=[\rm{diag}(h_j)]\in H$ this means
\[
fh=[h_1^{-1}f_1,\ldots,h_n^{-1}f_n].
\]
Also consider $J^\prime$ as a
principal $H$-bundle by identifying $H$ with $H_\fo$ via the epimorphism
\[
U_\fo\rightarrow  H, \quad f\mapsto  [\rm{diag}(f(O_j))],
\]
which has kernel $U_O$.  We will construct a principal bundle morphism 
$\mu:J^\prime\to\P Fr(E)$ such that $\mu(h^{-1}x)=\mu(x)h$
for $x\in J^\prime$ and $h\in H$. 
By (\ref{eq:E}), for each $L\in Pic(X)$ we have a natural map
$L^\times_\fo\to Fr(E)\vert L$ which assigns to a local section $\sigma$ its restrictions 
$(\sigma\vert O_1,\ldots,\sigma\vert O_{n+1})$. Now if $f\in U_O$ then $f\sigma$ is mapped to 
$(c\sigma\vert O_1,\ldots,c\sigma\vert O_{n+1})$ where $c$ is the value $f$ takes at every $O_j$. 
This induces a injection
\[
\mu: J^\prime\vert L\simeq L^\times_\fo/U_O\to \P Fr(E)\vert L\ ; 
\quad [\sigma]\mapsto [\sigma\vert O_1,\ldots,\sigma\vert O_{n+1}]
\]
where square brackets denote the respective equivalence classes. This map is clearly
$H$-equivariant. By globalising this we obtain an injective morphism of principal bundles. 

Now, $\pi^* Pic(\Xo)$ is trivialisable by its tautological section (and a little thought shows this
must actually be an algebraic section), whence $\P Fr(\pi^*E)$ is canonically trivialisable
and so are all the bundles $Gr_k(\pi^*E)$. $\Box$

It is an easy exercise to verify that, at each point $L^\prime\in Pic(X^\prime)$ with
$\pi(L^\prime)=L$, this canonical trivialisation maps
\[
\P\pi^*E\vert L^\prime\to \P^n,\quad [\sum v_js_j]\mapsto [v_1,\ldots,v_n],
\]
where $s_j\in L\vert O_j$ are such that $L^\prime$ is obtained from $L$ by identifying fibres so
that $s_j$ is identified with $s_1$.

For later use, let us now show that any section of either $\P Fr(\pi^*E)$ or
$Gr_k(\pi^*E)$ which is lifted from a section of $\P Fr(E)$ or $Gr_k(E)$ produces, 
by the trivialisation, an $H$-equivariant map of $J^\prime$ into $G^\C$ or $\Gr$. 
To show this we need to examine more closely the
canonical trivialisation. Let $\chi:J^\prime\to\pi^*J^\prime$ be the tautological section 
and lift $\mu$ to $\mu^\prime:\pi^*J^\prime\to\P Fr(\pi^*E)$. Then we have a section
\[
\tau =\mu^\prime\circ\chi:J^\prime\to\P Fr(\pi^*E)
\]
with the property $\tau(hx)=\tau(x)h^{-1}$ for $h\in H$, $x\in J^\prime$. This induces the canonical
trivialisation
\[
\P Fr(\pi^*E)\to J^\prime\times G^\C\ ; \quad \tau(x)g\mapsto (x,g).
\]
\begin{lem}
\label{lem:equivar}
Let $\sigma^\prime$ be a section of $\P Fr(\pi^*E)$ and $\psi:J^\prime\to G^\C$ be the map it
defines by $\sigma^\prime(x)=\tau(x)\psi(x)$. 
Then $\psi(hx)=h\psi(x)$ for all $h\in H$ if and only if
$\sigma^\prime$ is the lift of a section $\sigma$ of $\P Fr(E)$ over $J$.
\end{lem}
{\bf Proof.} Suppose $\sigma^\prime$ is such a lift, then $\sigma^\prime(hx)=\sigma^\prime$ 
for $h\in H$.  But
\[
\sigma^\prime(hx)=\tau(hx)\psi(hx)=\tau(x)h^{-1}\psi(hx)
\]
whence $h^{-1}\psi(hx)=\psi(x)$. Clearly the identities reverse to give us the converse. $\Box$

Notice that the lemma remains true if we replace $\P Fr(E)$ by $Gr_k(E)$ 
and $G^\C$ by $\Gr$. In this case we identify $Gr_k(\pi^*E)$ with $\P Fr(\pi^*E)\times_{G^\C}\Gr$ 
so that its canonical trivialisation can be written
\[
\P Fr(\pi^*E)\times_{G^\C}\Gr\to J^\prime\times \Gr\ ; \quad [\tau(x)g,V]\mapsto (x,gV) \ g\in 
G^\C.
\]
Again, square brackets denote an equivalence class. Thus a section $\sigma^\prime$ of $Gr_k(\pi^*E)$
defines a map $\psi:J^\prime\to\Gr$ by $\sigma^\prime(x)=[\tau(x),\psi(x)]$.

\subsection{Hermitian structure.}

So far we have not used the reality conditions on $(X,\lambda,\caL)$ at all: 
indeed, all we required were $n+1$
distinct points on $X$. The reality conditions arise because we must view $\Gr$ and $\PU$ as 
symmetric spaces. Each of these has, up to scaling, a unique group invariant metric. Because
the maps we construct will be {\it equi}-harmonic (i.e.\ harmonic with respect to any of these group
invariant metrics) the scaling will prove irrelevant. But the construction we are going to use, from
\cite{McI1}, works directly with $\Gamma(\caL)$ equipped with the inner product (\ref{eq:h}). So the aim
of this section is to show how this fits in with the fibres on $\pi^*E$ and the inner product they
obtain from the projective trivialisation.

First, we must restrict our attention to a particular real submanifold of $Pic(X^\prime)$:
\[
N^\prime = \{\caL^\prime\in Pic(X^\prime):\caL^\prime\otimes\overline{\rho_*\caL^\prime}\simeq
\caO_{X^\prime}(R)\}.
\]
Notice that $N^\prime$ consists of bundles of degree $g+n$  and is a translate in $Pic(X^\prime)$ of the
real subgroup
\[
J_R^\prime=\{L\in J^\prime:L^{\prime -1}\simeq\overline{\rho_*L^\prime}\}.
\]
Now for each $\caL^\prime\in N^\prime$ with $\caL=\pi(\caL^\prime)$ we have vector space isomorphisms
\begin{equation}
\label{eq:restriction}
\Gamma(\caL)\to \oplus_{j=1}^{n+1}\caL\vert O_j\to (\caL\vert O_1)^{n+1}
\end{equation}
where the first arrow is restriction to the fibres: this is bijective because
(we recall from section \ref{sec-spectral}) 
$\caL$ has no non-zero global sections vanishing at all the points
$O_j$ and $\dim\Gamma(\caL) = n+1$. The second arrow is the identification of fibres
$\caL\vert O_j\to\caL\vert O_1$ encoded in $\caL^\prime$. Indeed, this is how the canonical
trivialisation of $\P\pi^*E$ works. Notice that $\Gamma(\caL^\prime)$ is one dimensional since
it consists of all global sections of $\caL$ satisfying the $n$ fibre identifications. 
Therefore, given $\caL^\prime\in N^\prime$, we can equip 
$\Gamma(\caL)$ with another Hermitian inner product, unique up to positive real scaling, given by the formula
\begin{equation}
\label{eq:htilde}
\tilde h(\sigma_1,\sigma_2)=\sum_{j=1}^{n+1}(\frac{\sigma_1\vert O_j}{\sigma\vert O_j})
\overline{(\frac{\sigma_2\vert O_j}{\sigma\vert O_j})},
\end{equation}
where $\sigma\in\Gamma(\caL^\prime)$ is non-zero. It is this inner product which the canonical 
trivialisation of
$\pi^*E$ induces. The next lemma shows that provided we work with $\pi^*E\vert N^\prime$ the two
Hermitian structures agree. 
\begin{lem}
Let $\caL^\prime\in N^\prime$ and set $\caL=\pi(\caL)$. Then the two Hermitian inner products $h$ and
$\tilde h$ on $\Gamma(\caL)$, given by (\ref{eq:h}) and (\ref{eq:htilde}) respectively, agree (up to
real scaling).
\end{lem}
{\bf Proof.} Since $\caL^\prime$ has a non-zero global section, 
$\caL^\prime\simeq\caO_{X^\prime}(D)$ where 
$D$ is a positive divisor of degree $g+n$. This divisor $D$ is unique since $\dim\Gamma(\caL^\prime)=1$.
By definition of $N^\prime$, $D+\rho_*D-R$ is the divisor of a rational function $f$ taking the
same value at each point $O_j$. Indeed, we may choose $f$ so that $\overline{\rho_*f}=f$ and we 
may assume that $f$ is positive over $\vert\lambda\vert =1$ (\cite{McI2} lemma 2). Let
$\sigma\in\Gamma(\caL^\prime)\subset\Gamma(\caL)$ be a non-zero section, then $\sigma$ has divisor 
$D$ so on $\Gamma(\caL)$ we have
\[
h(\sigma_1,\sigma_2)  =  \sum_{j=1}^{n+1}f(O_j)
(\sigma_1\sigma^{-1})(O_j)\overline{\rho_*(\sigma_2\sigma^{-1})(O_j)}
   =  f(O_1)\tilde h(\sigma_1,\sigma_2), 
\]
and $f(O_1)$ is real and positive.$\Box$

From now on we will work exclusively with the bundle $E^\prime=\pi^*E\vert N^\prime$. Let
$\tau:N^\prime\to \P Fr(E^\prime)$ be its canonical projective framing and let $\P U(E^\prime)$
be the pullback of $N^\prime\times \PU$ by the trivialisation of $\P Fr(E^\prime)$ induced by $\tau$. 
Of course, this makes $\tau$ a section of $\P U(E^\prime)$.
This is a bundle of projective unitary groups and induces a Hermitian symmetric space structure on
each fibre of $Gr_k(E^\prime)$. From the
definition of $\tilde h$ and the previous lemma we see that this structure is projectively
equivalent with that induced by the fibrewise Hermitian metric $h$ on $E^\prime$. It follows
that the inclusion $J^\prime_R\hookrightarrow \P U(E^\prime)$ embeds $J^\prime_R\cap H_\fo$ into
$\PU$, hence $J^\prime_R$ is a real torus of dimension $g+n$.

\section{Harmonic maps into $\Gr$ and $\PU$.} 

For simplicity, set $G=\P U_{n+1}$.
We need to recall the characterisation of (equi)-harmonic maps into $\Gr$ and $\PU$ using maps into the
loop group $\Lambda G$ of real-analytic maps $S^1\to G$ (cf.\ \cite{BurFPP,Uhl}). 
Recall that $\Gr$ is isomorphic to the homogeneous space $G/K$ where $K$ is the fixed point
subgroup of the involution $g\mapsto \nu g\nu^{-1}$ determined by
\[
\nu = \left(\begin{array}{cc} I_k & 0 \\ 0 & -I_{n+1-k}\end {array}\right)
\]
where $I_k$ is the $k\times k$ identity matrix. 
This involution induces a symmetric decomposition of $\fg$,
the Lie algebra of $G$, into $\fg=\fk +\fm$, where $\fk$ is the Lie algebra for $K$. Recall (from e.g.\
\cite{BurP}) that a (based) map 
\[
\varphi:\R^2\to G/K,\quad f(0)=1\cdot K\]
is harmonic iff it possesses an extended frame $F_\zeta:\R^2\to \Lambda 
G$ (i.e.\ $F_1$ frames $\varphi$) for which 
\begin{equation}
\label{eq:Fz}
F^{-1}_\zeta\partial F_\zeta/\partial z = \zeta^{-1}A_{-1} +A_0, 
\end{equation}
where $z=x+iy$, $A_0$ takes values in $\fk^\C$($=\fk\otimes\C$) and $A_{-1}$ takes values in $\fm^\C$. 
Now let $\fm_-^\C$ and $\fm_+^\C$ denote the subspaces of, respectively, lower triangular and upper triangular 
matrices in $\fm^\C$. Then a useful corollary to the above fact is:
\begin{lem}
\label{lem:Grharmonic}
Let $F_\lambda:\R^2\to\Lambda G$ satisfy
\[
F_\lambda^{-1}\partial F_\lambda/\partial z = \lambda^{-1}a_{-1} +a_0
\]
where $a_0$ takes values in $\fk^\C+\fm_+^\C$ and $a_{-1}$ takes values in $\fm_-^\C$. Then $F_1$ frames a
harmonic map into $\Gr$.
\end{lem}
{\bf Proof.} Define 
\[
\kappa = \left(\begin{array}{cc} I_k & 0 \\ 0 & \zeta I_{n+1-k}\end{array}\right),
\]
where $\zeta=\sqrt\lambda$.  One easily checks that $\tilde F_\zeta=\kappa^{-1}F_{\zeta^2}\kappa$
satisfies (\ref{eq:Fz}) and therefore $F_1=\tilde F_1$ frames a harmonic map into $\Gr$.$\Box$

We will also need the analogous characterisation of harmonic maps $\varphi:\R^2\to G$. Recall from
\cite{Uhl} that $\varphi$ is harmonic if and only if it
possesses an extended frame $F_\lambda:\R^2\to\Lambda G$ for which
\begin{equation}
\label{eq:Gharmonic}
F^{-1}_\lambda\partial F_\lambda/\partial z = (1-\lambda^{-1})A
\end{equation}
where $A$ takes values in $\fg^\C$.

\subsection{Maps into $\Gr$.}

As explained in the introduction, when we choose
the right section of $Gr(E^\prime)$ restricted to the right directions along $J^\prime$ we
will obtain harmonic maps. These choices are made as follows.

We assume, without loss of generality, that $k\leq (n+1)/2$.
Let $(X,\lambda,\caL)$ be spectral data as described as above
and let us further insist that $\lambda$ has
at least $k$ zeroes $P_1,\ldots,P_k$ of order two and therefore $k$ poles $Q_1,\ldots,Q_k$ of order
two given by $Q_j=\rho(P_j)$. Define positive divisors $D_0,D_\infty$ of degree $n+1-k$ by
\[
D_\infty=(\lambda)_\infty-Q_1-\ldots -Q_k, 
\]
where $(\lambda)_\infty$ is the divisor of poles
of $\lambda$, and
$D_0=\rho_*D_\infty$. By earlier remarks, $\caL(-(\lambda)_\infty)$ is non-special and it follows by
Riemann-Roch that $\Gamma(\caL(-D_\infty))$ has dimension $k$. So $D_\infty$ determines a rank $k$
subbundle of the vector bundle $E^\prime$ and therefore a section of $Gr_k(E^\prime)$. By applying the
canonical trivialisation we obtain a map $N^\prime\to\Gr$. For each $\caL^\prime\in N^\prime$
we have an isomorphism $J^\prime_R\to N^\prime$ and the combination of these two maps we will
define to be
\[
\psi_{\caL^\prime}:J^\prime_R\rightarrow \Gr.
\] 
Notice that the spectral data only determines this map up to the choice of $\caL^\prime$ over $\caL$. 
However, the group $H_\fo\cap J_R^\prime$  acts transitively on each fibre of $\pi:N^\prime\to N$ and this
action induces isometries after canonical trivialisation of $\P Fr(E^\prime)$, 
so the map $\psi_{\caL^\prime}$ is determined by the spectral data up to isometries of the target. 

Now we will define the real homomorphism $\gamma:\R^{2k}\to J^\prime_R$ for which $\varphi=\psi_{\caL^\prime}
\circ\gamma$
is pluri-harmonic. Set $X^\prime_o= X^\prime-\{O\}$ and recall 
(from e.g.\ \cite{Ser}) that there is, for any point $B\in X^\prime_o$, an Abel map 
\[
\caA_{B}^\prime:X^\prime_o\rightarrow J^\prime \simeq \Gamma(\Omega_{X^\prime})^*
/ H_1(X^\prime_o,\Z) 
,\qquad P\mapsto \int_{B}^{P}
\]
with the integral considered modulo periods. About each $P_j$ we have a local
coordinate $\zeta_j=\sqrt\lambda$. We characterise $\gamma$ by
\begin{equation}
\label{eq:tangent}
\partial\gamma/\partial z_j\vert_{z=0}=\partial\caA^\prime_{P_j}/\partial\zeta_j\vert_{\zeta_j=0},
\end{equation}
where $z=(z_1,\ldots,z_k)$ are the standard complex coordinates on $\R^{2k}$. 
\begin{thm}
\label{th:Gr}
For each $\caL^\prime\in N^\prime$ the map $\varphi=\psi_{\caL^\prime}
\circ\gamma:\R^{2k}\to\Gr$ is pluri-harmonic. Moreover, $\varphi(z_m)$ 
is non-conformal for each $m=1,\ldots,k$.
\end{thm}
Before presenting the proof we need a few preliminary results.
\begin{lem}
\label{lem:immersion}
The real homomorphism $\gamma:\R^{2k}\to J^\prime_R$ is an immersion.
\end{lem}
{\bf Proof.} Notice that $\gamma$ is tangent at the identity to 
$\caA^\prime_{D_0}:X_o^{\prime(n+1-k)}\to J^\prime$, the Abel map on the symmetric product of
$X_o^\prime$. It suffices to show that $d\caA^\prime_{D_0}$ is an injection at the 
base point $D_0$.
The fibre of $\caA_{D_0}^\prime$ over the identity is the
projective space $\P V$ where
\[
V=\{f\in \C(X): (f)\geq -D_0\ \hbox{and $f(O_j)=f(O_1)\neq 0$ for every
$j$}\}.
\]
Any $f\in V$ has degree $\leq k$, so unless it is constant it cannot
take the same value at the $n+1$ distinct points $O_j$. Therefore $V=\C$ and it follows that
$d\caA^\prime_{D_0}$ is injective at $D_0$.$\Box$

Next we must describe the image of $\gamma$ as a family of line bundles. Set 
$X^\prime_A=X^\prime-\lambda^{-1}(\{0,\infty\})$ and for each $j=1,\ldots k$
let $U_j$ be an open coordinate disc about $P_j$ for the coordinate $\zeta_j$. Further, let $V$ be a union of
coordinate discs about the other points at $\lambda=0$ and assume $U_1,\ldots,U_k,V$ are mutually disjoint. Then 
\[
\caU = \{ X^\prime_A,U_1,\ldots,U_k,V,\rho_*U_1,\ldots,\rho_*U_k,\rho_*V\}
\]
is a Leray cover for $X^\prime$ and any line bundle over it can be prescribed by a 1-cocycle for
the open cover $\caU$.
\begin{lem}
Let $[c]:\R^{2k}\to H^1(\caO_{X^\prime}^\times)$ be the map which gives the class of the 1-cocycle defined by
\[
c(z)=\left\{\begin{array}{ll} \exp(z_j\zeta_j^{-1}) & \hbox{on}\ X^\prime_A\cap U_j,\ j=1,\ldots,k \\
                              1 & \hbox{on}\ X^\prime_A\cap V \end{array}\right.
\]
and $\overline{\rho_*c(z)}=c(z)^{-1}$. Then $[c(z)]$ is the class of $\gamma(z)$.
\end{lem}
{\bf Proof.} Clearly $[c(z)]$ is the class of a degree zero line bundle in $J^\prime_R$. Consider, for each
$m=1,\ldots,k$, 
\[
\partial [c]/\partial z_m\vert_{z=0}=
\left\{\begin{array}{ll} \zeta_m^{-1} & \hbox{on}\ X^\prime_A\cap U_m\\
                              0 & \hbox{on every other 2-simplex from}\ \caU.  \end{array}\right.
\]
Recall (e.g.\ from \cite{Ser}) that the Serre duality $H^1(\caO_{X^\prime})\simeq \Gamma(\Omega_{X^\prime})^*$
identifies the class of this 1-cocycle with the residue 
map $f:\omega\mapsto res_{P_m}\zeta^{-1}_m\omega$ for all
$\omega\in\Gamma(\Omega_{X^\prime})$. But
\[
res_{P_m}\zeta_m^{-1}\omega = (\omega/d\zeta_m)(P_m) =( \frac{\partial}{\partial
\zeta_m}\int_{P_m}^{\zeta_m}\omega)\vert_{\zeta_m=0}.
\]
It follows that 
$(\partial [c]/\partial z_m)(0) = (\partial\caA_{P_m}/\partial \zeta_m)(0)$
whence $[c]$ and $\gamma$ are the same map.$\Box$

Finally, before we prove theorem \ref{th:Gr} let us make a few crucial observations concerning
the canonical trivialisation of $Gr_k(E^\prime)$ over $\caL^\prime\otimes\gamma(z)$.

Any line bundle $\caL^\prime$ over $X^\prime$ can be thought of as a line bundle $\caL$ over $X$ equipped with
fibre identifications $g_j(z):L(z)\vert O_{j}\to L(z)\vert O_1$, or equally, non-zero elements
$s_j\in\caL\vert O_j$, only fixed up to common scaling, for which $g_j$ identifies $s_j$ with $s_1$. From this
point of view $\gamma(z)$ is the line bundle $L(z)=\pi(\gamma(z))$ with non-zero elements $s_j(z)\in L(z)\vert
O_j$ which can be determined as follows. The 1-cocycle $c(z)$ equips $\gamma(z)$ with 
trivialising sections over each open set in
$\caU$. Let $s_z$ be the section over $X_A^\prime$: it is only determined up to scaling. Then $s_z$ is
a section of $L(z)$ over $X_A$ and $s_j(z) = s_z\vert O_j$.  
Notice that we may assume that $\overline{\rho_*s_z}=s_z^{-1}$. 
Now for any $\caL^\prime\in N^\prime$, with non-zero global section $\sigma$ (recall this is unique up to
scaling) the canonical trivialisation of $Gr_k(E^\prime)$ over $\caL^\prime\otimes\gamma(z)$ works as follows.
Set $\caL=\pi(\caL^\prime)$ and $\caL_z=\caL\otimes L(z)$, then
\begin{equation}
\label{canon1}
Gr_k(\Gamma(\caL_z))\rightarrow \Gr;\quad e_1\wedge\ldots\wedge e_k\mapsto v_1\wedge\ldots\wedge v_k,
\end{equation}
where $v_j\in\C^{n+1}$ is given by
\begin{equation}
\label{eq:canon2}
v_j =((e_j/(\sigma\otimes s_z))(O_1),\ldots,(e_j/(\sigma\otimes s_z))(O_{n+1}))
\end{equation} 

\noindent
{\bf Proof of Theorem \ref{th:Gr}.} The aim is to produce a map $F_\lambda:\R^{2k}\to \Lambda G$ 
for which $F_1$ frames $\varphi$
and satisfies (\ref{lem:Grharmonic}) for each $z_m$, $m=1,\ldots,k$. 
Set $\caL_z=\caL\otimes L(z)$ for each $z\in\R^{2k}$.  Let
$e_1^z,\ldots,e_{n+1}^z$ be a unitary frame for $\Gamma(\caL_z)$ for which $e_1^z,\ldots,e_k^z$ span the
$k$-plane $V=\Gamma(\caL_z(-D_\infty))$ and $e_{k+1}^z,\ldots,e_{n+1}^z$ span $V^\perp=
\Gamma(\caL_z(-\sum_{j=1}^k
P_j))$. A straightforward calculation with the trace norm (see [10] lemma 6) shows that $e,e^\prime\in
\Gamma(\caL_z)$ are orthogonal whenever $e^\prime\otimes\overline{\rho_*e}$ has divisor of zeroes at
least $(\lambda)_0$ (or $(\lambda)_\infty$), therefore these two
spaces are orthogonal.
Since any bundle over $\R^{2k}$ is analytically trivial we can choose the $e_j^z$ to be analytic in $z$.
Then $e_1^z\wedge\ldots\wedge e_k^z$ is the section of $Gr_k(E^\prime)$ (pulled back to $\R^{2k}$) which
produces $\psi_{\caL^\prime}$. 
 
Now let $f^z_1,\ldots,f^z_{n+1}$ be the
dual basis for $\Gamma(\caL_z)^*$, and define the matrix $F_\lambda(z)$ to have entries
\begin{equation}
\label{eq:frame}
F_\lambda(z)_{ij} = f^0_i e^z_j s_z^{-1}.
\end{equation}
To make sense of this, let $B$ denote the ring $\C[\lambda,\lambda^{-1}]$ and think of 
the $e_j^z$ as generators of the $B$-module $M_z=
\Gamma(X_A,\caL_z)$. This is free of rank $n+1$ since $\lambda_*\caL_z$ is a trivial 
rank $n+1$ bundle over $\P^1$.
We think of $f^0_j$ as generators for the dual
module so that the entries of $F_\lambda(z)$ are elements of $B$. By the reality
conditions $F_\lambda:\R^{2k}\to U_{n+1}$ for $\vert\lambda\vert =1$. I claim that
the projective frame $[F_1(z)]$ 
frames $\varphi(z)$, after possibly an isometry of $\Gr$. To see this notice that the only
difference between (\ref{eq:frame}) at $\lambda =1$ and (\ref{eq:canon2}) is that the former identifies
$\oplus\caL\vert O_j$ with $\C^{n+1}$ using the dual basis $\{f^0_j\}$ while the latter does this using the
fibre elements $\{\sigma\vert O_j\}$. Either identification respects the Hermitian structure (up to
scaling), hence the difference between them is a projective unitary transformation.  

Now observe that  
\begin{equation}
\label{eq:derivative}
(F_\lambda^{-1}\partial F_\lambda/\partial z_m)_{ij} = 
s_zf_i^z\partial(e_j^zs^{-1}_z)/ \partial z_m.
\end{equation}
If $e^z\in \Gamma(\caL_z(-D))$ is any smooth family
of sections, for some divisor $D$, define 
\begin{equation}
\label{eq:covariant}
\nabla_{z_m}e^z = 
s_z\partial(e^zs^{-1}_z)/ \partial z_m .
\end{equation}
Let $\pi_V^\perp$ be $B$-linear projection onto
$V^\perp\otimes B$ with kernel $V\otimes B$. Then $F_\lambda$ will 
satisfy the
conditions of lemma \ref{lem:Grharmonic} for each $z=z_m$, $m=1,\ldots,k$, if
\begin{equation}
\label{eq:perp}
\pi_V^\perp(\nabla_{z_m}e) \in\left\{\begin{array}{cc} 
\lambda^{-1}V^\perp&
\hbox{for}\ e\in V = \Gamma(\caL_z(-D_\infty)); \\ 
\Gamma(\caL_z) &\hbox{for}\ e\in V^\perp= \Gamma(\caL_z(-\sum P_j)).
\end{array}\right.
\end{equation}
To see that this holds we recall from \cite{McI1} lemma 7 that, whenever $e\in\Gamma(\caL_z(-D))$ for
some divisor $D$, then $\nabla_{z_m}e\in\Gamma(\caL_z(P_m-D))$. So for $e\in V$ 
\[
\nabla_{z_m}e\in\Gamma(\caL_z(P_m-D_\infty))=\C\langle \lambda^{-1}e_0,e_1,\ldots,e_k\rangle
\]
where $e_0\in V^\perp$ has divisor of zeros $(\lambda)_0-P_m$ (which is positive since $P_m$ is a 
ramification point of $\lambda$). Therefore $\pi_V^\perp(\nabla_{z_m} e)=\lambda^{-1}e_0$. On the other
hand, for $e\in V^\perp$ we have 
\[
\nabla_{z_m} e\in\Gamma(\caL_z(-\sum_{j\neq m} P_j))\subset\Gamma(\caL_z).
\]
Therefore (\ref{eq:perp}) holds, so  by lemma \ref{lem:Grharmonic}
$\varphi$ is harmonic is each variable $z_m$. 

To see that $\varphi(z_m)$ is non-conformal we will show that the map 
\[
A_m=\pi_V\nabla_{\bar z_m}\pi_V^\perp\nabla_{\bar z_m}:V\to V 
\]
has nowhere vanishing trace. This shows that
$(\partial\varphi/\partial\bar z_m,\partial\varphi/\partial\bar z_m)^\C$ vanishes nowhere, 
where this is the $(0,2)$-component of the
complexification of the metric induced on $\R^2$ by $\varphi$. First observe that $\nabla_{\bar
z_m}V\subset\Gamma(\caL_z(Q_m-D_\infty))$ so that we have
\begin{equation}
\label{eq:V}
V= \ker(\pi^\perp_V\nabla_{\bar z_m}) \oplus\C\langle \sigma^z_0\rangle
\end{equation}
where $\sigma_0^z$ generates $\Gamma(\caL_z(-\hat Q_m))$ and $\hat Q_m$ is 
the positive degree $n$ divisor $(\lambda)_\infty-Q_m$. Now observe that
\[
\nabla_{\bar z_m}\sigma_0^z\in V\oplus \C\langle\sigma_1\rangle
\]
where $\sigma_1$ generates $\Gamma(\caL_z(-D_m-\sum P_j))\subset V^\perp$ 
and $D_m$ is a degree $k-1$ positive 
divisor $D_\infty -Q_m$ since $Q_m$ is a ramification point of $\lambda$. 
Therefore we can write
\[
\pi_V^\perp\nabla_{\bar z_m}\sigma_0^z = a\sigma_1^z
\]
where $a$ depends only on $z$.  Next, observe that
\begin{equation}
\label{eq:nabla1}
\nabla_{\bar z_m}\sigma_1^z\in \Gamma(\caL_z(Q_m-D_m-\sum P_j)) = 
V^\perp\oplus\C\langle \lambda\sigma_0^z\rangle.
\end{equation}
It can be shown, by a computation almost identical to the one in \cite{McI1}, pp 849-850, that 
\[
\pi_V\nabla_{\bar z_m}\sigma_1^z = a^{-1}\lambda\sigma_0^z.
\]
Finally, when we pass from $V\otimes B$ to $V$ (by evaluation at $\lambda =1$) we see,
by choosing a basis for $V$ compatible with the splitting in (\ref{eq:V}), that
the trace of $A_m:V\to V$ is 1. $\Box$

\smallskip\noindent
{\em Remarks.} 1. In fact similar methods can be used to show that
\[
\Tr( \pi_V\nabla_{\bar z_l}\pi_V^\perp\nabla_{\bar z_m}) =0\ \hbox{for}\ l\neq m,
\]
and therefore if we restrict $\varphi$ to the 2-plane with tangent
$\sum a_m\partial/\partial z_m$ then $\varphi$ is non-conformal provided $\sum a_m^2\neq 0$. 

\smallskip\noindent
2. What happens if $P_m$ has ramification index greater than 1? In that case the divisor $D_m-Q_m$
in (\ref{eq:nabla1}) is positive, so that (\ref{eq:nabla1}) becomes $\nabla_{\bar
z_m}\sigma_1^z\subset V^\perp$. Hence $\Tr(\pi_V\nabla_{\bar z_m}\pi^\perp_V\nabla_{\bar z_m})=0$
i.e.\ $\varphi(z_m)$ is a minimal surface (i.e.\ conformal harmonic). If every 
$P_1,\ldots,P_k$ has this property then $\varphi$ is minimal when restricted to any 2-plane. Indeed,
we can raise the isotropy order (in the sense of \cite{Wood}) of $\varphi(z_m)$ further by increasing the ramification index of
$P_m$, in the same way we see this in \cite{McI1} Thm 4, p847.

\subsection{Maps into $\PU$.}

In this case for the spectral data $(X,\lambda)$ we do not assume that $\lambda$ has any 
double zeroes. Indeed, for simplicity we will assume that $\lambda$ only has simple zeroes
$P_1,\ldots,P_{n+1}$ with corresponding poles $Q_j=\rho_(P_j)$. But the singularisation $X^\prime$ we
will use is not the same as in the previous section. We want $X^\prime$ to have {\it two} singular
points $O$ and $S$, corresponding to two divisors: $\fo=O_1+\ldots+O_{n+1}$, the divisor of
zeroes of $\lambda-1$, and $\fs=S_1+\ldots+S_{n+1}$, the divisor of zeroes of $\lambda +1$. Each
divisor provides us with a singularisation of X: call them $X_\fo$ and $X_\fs$ respectively. Then
$X^\prime$ is a further singularisation of each of these. It follows that we have a commuting
diagram of generalised Jacobi varieties:
\[
\begin{array}{rcccl}
       &        & J^\prime &       &          \\
       &\swarrow &            &\searrow&         \\       
J_\fo  &         &\ \downarrow \pi          &          & J_\fs\\
       &\pi_\fo\searrow &            &\swarrow\pi_\fs&         \\       
       &        &      J    &       &          
\end{array}
\]
Both $\pi_\fo^*E$ and $\pi_\fs^*E$ are projectively trivial. Let us
denote by $\tau_O$ and $\tau_S$ the two sections of $\P Fr(\pi^*E)$ obtained by pulling back,
respectively, the canonical sections of $\P Fr(\pi_\fo^*E)$ and $\P Fr(\pi^*_\fs E)$. If we now
restrict our attention to $E^\prime=\pi^*E\vert N^\prime$ then $\tau_O$ and $\tau_S$ are sections
of $\P U(E^\prime)$. It follows that we have a change of frame map
\[
\Psi:N^\prime \to \PU, \quad \tau_O=\Psi\tau_S.
\]
For each $\caL^\prime\in N^\prime$ we define $\Psi_{\caL^\prime}$ by precomposing $\Psi$ with the
isomorphism $J_R^\prime \simeq N^\prime$ which identifies the identity in $J^\prime_R$ with
$\caL^\prime$.  

Since each zero $P_j$ of $\lambda$ is assumed to be simple we can use $\lambda$ as a local
coordinate (which we name $\zeta_j$) about each $P_j$. Define $\gamma:\R^{2n}\to J^\prime_R$ to
be the unique homomorphism satisfying (\ref{eq:tangent}) for $j=1,\ldots,n$.
Then $\gamma$ is an immersion by lemma \ref{lem:immersion} (adapted to the case where $X^\prime$
has two singular points).
\begin{thm}
\label{th:PU}
For each $\caL^\prime\in N^\prime$, the map 
$\varphi=\Psi_{\caL^\prime}\circ\gamma:\R^{2n}\to \PU$ is
pluri-harmonic. Moreover, $\varphi(z_m)$ is non-conformal for $m=1,\ldots,n$.
\end{thm}
{\bf Proof.} The proof is much the same as the proof of theorem \ref{th:Gr}. Let $L:\R^{2n}\to J$
be defined by $L=\pi\circ\gamma$, set $\caL=\pi(\caL^\prime)$ and define $\caL_z=\caL\otimes L(z)$. 
Let $\{e_i^z\}$ be a unitary basis for $\Gamma(\caL_z)$ 
and take $\{f_j^z\}\subset\Gamma(\caL_z)^*$ to be the dual basis. As earlier, $L(z)$ admits 
admits a trivialising section $s_z$
over $X-\lambda^{-1}(\{0,\infty\})$ with the properties: (i) $s_z\exp(z_j\lambda^{-1})$ is
non-vanishing and holomorphic about $P_j$ while $s_z\exp(-\bar
z_j\lambda)$ is holomorphic and non-vanishing about $Q_j$, (ii)
$\overline{\rho_*s_z}=s_z^{-1}$. 
As before we introduce the matrix  $F_\lambda(z)$ with entries
\[
F_\lambda(z)_{ij} = (f^0_ie_j^zs_z^{-1}).
\]
I claim that: (a) after possibly an isometry, $\varphi$ is framed in $\PU$ 
by $[F_{-1} F_1^{-1}]$, where the square brackets denote the class in $\PU$;
(b) $[F_\lambda F_1^{-1}]$ is an extended frame for a
pluri-harmonic map of $\R^{2n}$ into $\PU$.

For (a), let $[F_\lambda]$ denote the map into $\PU$ at each $\vert
\lambda\vert = 1$ and let $[e_z]$ denote the projective frame determined by the
$e_j^z$. By the same reasoning as in the proof of the previous theorem,
there exist $a,b\in\PU$ for which
\[
[e_z]=\tau_O(\gamma(z))a[F_1(z)]= \tau_S(\gamma(z))b[F_{-1}(z)]
\]
for all $z$, so that $\Psi\circ\gamma = b[F_{-1}F_1]a^{-1}$. Since the metric on $\PU$ is bi-invariant
$[F_{-1}F_1^{-1}]$ differs from $\varphi$ by an isometry.

For (b), we compute, for each $m=1,\ldots,n$,
\[
F_\lambda^{-1}\partial F_\lambda/\partial z_m = (a^{(m)}_{ij})\ ; \quad 
a^{(m)}_{ij} =
f^z_i\nabla_{z_m}e_j^z
\]
(cf.\ (\ref{eq:derivative}) and (\ref{eq:covariant})). A computation similar 
to the one earlier shows that
\begin{equation}
\label{eq:nabla2}
\nabla_{z_m}e \in\left\{\begin{array}{cc} \Gamma(\caL_z) &
\hbox{for}\ e\in \Gamma(\caL_z(-P_m)); \\ \lambda^{-1}\Gamma(\caL_z)\oplus\Gamma(\caL_z) &
\hbox{otherwise.}
\end{array}\right.
\end{equation}
It follows that
$F_\lambda^{-1}\partial F_\lambda/\partial z_m$ has linear dependence on
$\lambda^{-1}$. So if we define $\tilde F_\lambda = F_\lambda F_1^{-1}$ then
\[
\tilde F_\lambda^{-1}\frac{\partial \tilde F_\lambda}{\partial z_m}
 = (1-\lambda^{-1})\tilde A_m 
\]
for some matrix function $\tilde A_m$ independent of $\lambda$, whence the map is harmonic. 

Now we will check that $\varphi(z_m)$ is non-conformal. Notice that means checking that
\[
(\partial\varphi/\partial z_m,\partial\varphi/\partial z_m)^\C = \Tr(\tilde A_m^2)
\]
is not identically zero. First, let us note that if we write
\[
F^{-1}_\lambda\partial F_\lambda/\partial z_m = \lambda^{-1} A_m + B_m
\]
then a simple calculation shows that $\tilde A_m = \rm{Ad}F_1(-2A_m)$ so that $\Tr(\tilde
A_m^2)=\Tr(4A_m^2)$. Now fix the basis $\{e_j^z\}$ so that $e_j^z\in 
V=\Gamma(\caL_z(-P_m))$ for $j\neq
m$ and $e_m^j$ generates $V^\perp=\Gamma(\caL_z(-\hat Q_m))$, where $\hat Q_m$ is the positive
divisor $(\lambda)_\infty-Q_m$. Notice from (\ref{eq:nabla2}) that 
$\nabla_{z_m}V\subset \Gamma(\caL_z)$. This means that, 
in this basis, all the columns of $A_m$ except possibly the $m$-th are zero. Therefore,
$\Tr(A_m^2)=A_{m,mm}^2$, where $A_{m,mm}$ is the $m$-th entry of the $m$-th column.
To compute this, we first notice that 
\[
\nabla_{z_m}e^z_m\in\Gamma(\caL_z(P_m-\hat Q_m))=\C\langle\lambda^{-1}e^z,e^z_m\rangle.
\]
where $e^z$ generates $\Gamma(\caL_z(-\hat P_m))$ with $\hat P_m=(\lambda)_0-P_m$. But 
\begin{equation}
\label{eq:Pm}
\Gamma(\caL_z(-\hat P_m))\not\subset V
\end{equation}
so that $f^z_m(e^z)\neq 0$, therefore $A_{m,mm}\neq 0$ whence $\varphi(z_m)$ is non-conformal.
$\Box$

\smallskip\noindent
{\em Remark.} One knows from \cite{BurFPP} that all non-conformal tori in
$\mathbb{P}U_2$ are of semisimple finite type and I claim the construction
above will give them all. But in general it is not clear to what
extent the maps of semisimple finite type account for the harmonic tori.
As with the case of maps into $\Gr$, we can increase the isotropy order of $\varphi(z_m)$ by
increasing the ramification index of $P_m$. In particular, $\varphi(z_m)$ is minimal whenever $P_m$
is a ramification point. For in that case $\hat P_m$ has support at $P_m$ so (\ref{eq:Pm}) 
becomes $\Gamma(\caL_z(-\hat P_m))\subset V= \ker f^z_m$, whence $A_{m,mm} = 0$.

\subsection{Equivariant maps.}

For any subgroup $S\subset\R^{2k}$ we will say a map $\varphi:\R^{2k}\to \Gr$ is $S$-equivariant if
there exists a homomorphism $h:S\to \PU$ for which $\varphi(z+s)=h(s)\varphi(z)$ for all
$z\in \R^{2k}$ and $s\in S$. 
\begin{prop}
\label{pp:equiv}
Let $\varphi=\psi\circ\gamma:\R^{2k}\to J^\prime_R\to\Gr$ be given by theorem \ref{th:Gr} and let
$S\subset\R^{2k}$ be the subgroup covering $\gamma(\R^{2k})\cap H_\fo$.
Then $\varphi$ is $S$-equivariant.
\end{prop}
{\bf Proof.} This follows at once from lemma \ref{lem:equivar}. From the proof of theorem
\ref{th:Gr} we know that $\psi$ is derived by canonical trivialisation of a section of $Gr_k(E)$
lifted to $Gr_k(E^\prime)$. By lemma \ref{lem:equivar} and the definition of $S$, 
$\varphi(z+s)=\gamma(s)\psi(\gamma(z))$ for $z\in\R^{2k}$ and $s\in S$, 
where we have assumed the identification of $H_\fo$ with
$H\subset\PU$.$\Box$

Clearly the same result holds for the maps $\varphi:\R^{2n}\to\PU$ constructed by theorem 
\ref{th:PU}, where now we have a homomorphism $S\to {\rm Isom}(\PU)\simeq \PU\times\PU$. 
We will say the map $\varphi$ is totally equivariant if it has maximal equivariance group
($\R^{2k}$ or $\R^{2n}$ as appropriate).
\begin{cor}
\label{cor:equiv}
(i) If the spectral curve $X$ for $\varphi$ has genus $g< 2k$ (resp.\ $g< 2n$) then $\varphi$ is
an $\R^{2k-g}$-equivariant map into $\Gr$ (resp.\ an $\R^{2n-g}$-equivariant map into $\PU$).
In particular, when $X$ is the Riemann sphere the map is totally
equivariant. 
\end{cor}
{\bf Proof.} The kernel of $d(\pi\circ\gamma)$ has dimension $2k-g$ (resp.\ $2n-g$).$\Box$

\section{The $\theta$-line bundle and algebraic maps.}

The aim of this section is to show that the maps produced from the 
canonical trivialisation can be written down using sections of the pullback of the $\theta$-line
bundle to $J^\prime$. We then look at the possibility that the $\varphi$ is algebraic (which happens
precisely when $\gamma$ is algebraic). It turns out, at least for maps into $\CP^n$, that this forces
$\varphi$ to be doubly periodic and therefore give harmonic 2-tori.

\subsection{The pullback of the $\theta$-line bundle to $J^\prime$.}

Let $\caT$ denote the $\theta$-line bundle over $J$. It has, up to scaling, a unique global section
$\theta$ which we will also think of as the classical Riemann $\theta$-function over $\C^g$ i.e.\ we
will assume we have made all the appropriate identifications so that $J\simeq\C^g/\Lambda$ where
$\Lambda\simeq H_1(X,\Z)$. We want to describe $\Gamma(\pi^*\caT)$ explicitly. General theory (see
e.g.\ \cite{Har} p 128) tells us that the space $\Gamma(\pi^*\caT)$ of globally {\em algebraic} sections
should be isomorphic to the coordinate
ring of any fibre of $\pi:J^\prime\to J$ as a vector space. We will actually compute a spanning set 
of sections for $\Gamma(\pi^*\caT)$.

First, notice that we can describe $J^\prime$ as a fibre product over $J$:
\[
J^\prime\simeq J_1\times_J \ldots \times_J J_n,
\]
where $J_k$ is the generalised Jacobian for the modulus $O_1+O_{k+1}$. Each $J_k$ is a
$\C^\times$-bundle over $J$ once we have fixed an isomorphism of $H_k\simeq J_k/J$ with $\C^\times$. 
We fix this so that for any $L\in J$ the group 
$\C^\times$ acts on ${\rm Hom}(L\vert O_{k+1},L\vert O_1)$ by
multiplication. Every $\C^\times$-bundle is characterised by its pullback to $X$ via the Abel
map $\caA:X\to J$ and with this choice $J_k$ is characterised by 
$\caA^*J_k=\caA(O_1-O_{k+1})$ (cf.\ \cite{Ser} VII no 9). For each $k$ let 
$\pi_k:J^\prime\to J_k$ denote the natural fibration. To describe
$\Gamma(\pi^*\caT)$ we need to introduce $t_k:J\to J$, the translation by $\caA(O_1-O_{k+1})$, and
$s_k$, the tautological section of $\pi^*J_k\to J^\prime$. We will denote by $t_k^m$ the $m$-th
iterate of the translation $t_k$ (i.e.\ translation by $\caA(mO_1-mO_{k+1})$). 
\begin{prop}
There is a vector space isomorphism between $\C[x_1,x_1^{-1}]\otimes\ldots\otimes\C[x_n,x_n^{-1}]$ 
and $\Gamma(\pi^*\caT)$ given by
\begin{equation}
\label{eq:Theta}
x_1^{k_1}\ldots x_n^{k_n}\mapsto \theta_{{k_1},\ldots,{k_n}} = s_1^{k_1}\otimes\ldots\otimes
s_n^{k_n}\otimes\pi^*(t_1^{-k_1}\ldots t_n^{-k_n})^*\theta .
\end{equation}
\end{prop}
{\bf Proof.} First we must establish that every $\theta_{{k_1},\ldots,{k_n}}$ belongs to
$\Gamma(\pi^*\caT)$. To simplify notation, define $\xi = \sum_{j=1}^n\caA(k_jO_{j+1}-k_jO_1)$, then 
$t_\xi=t_1^{-k_1}\ldots t_n^{-k_n}$. Let us show that
\begin{equation}
\label{eq:caT}
t_\xi^*\caT \simeq\caT\otimes J_1^{k_1}\otimes\ldots\otimes J_n^{k_n},
\end{equation}
where, with a slight abuse of notation, we are thinking of each $J_m\to J$ as the associated
line bundle to the $\C^\times$-bundle.
Since the pullback by the Abel map induces an isomorphism 
$H^1(X,\caO_X^\times)\simeq H^1(J,\caO_J^\times)$ 
it suffices to show that 
$\caA^*(t_\xi^*\caT)=\caA^*\caT\otimes\xi$, where we use $\xi$ here to denote the line bundle over
$X$ for the point $\xi\in J$. But this follows immediately from Riemann's vanishing theorem.

Now if we pull back (\ref{eq:caT}) to $J^\prime$ and use the fact that $s_k$ trivialises $\pi^*J_k$
we see that $\pi^*t_\xi^*\caT\simeq\caT$ and tensoring by $s_1^{k_1}\otimes\ldots\otimes
s_n^{k_n}$ identifies $\pi^*t_\xi^*\theta$ with $\theta_{{k_1},\ldots,{k_n}}$.

To prove that these span $\Gamma(\pi^*\caT)$ we take the point of view that each section is a
holomorphic function on the universal cover $\C^{g+n}$ of $J^\prime$ satisfying certain
functional equations. To identify $J^\prime$ with $\C^{g+n}/\Lambda^\prime$, where 
$\Lambda^\prime\simeq H^1(X^\prime-\{O\},\Z)$, we augment the standard first homology 
basis $\{a_j,b_j:j=1,\ldots,g\}$ of $X$
with a positively oriented cycle $a_{g+j}$ enclosing a small disc about each 
$O_{j+1}$, $j=1,\ldots,n$. Then $\{\oint_{a_j}:j=1,\ldots,g+n\}$ is a 
basis for $\Gamma(\Omega_{X^\prime})^*$. A
function $f:\C^{g+n}\to\C$ corresponds to a section of $\pi^*\caT$ if and only if
\begin{equation}
\label{eq:functional}
(f/\pi^*\theta)(\tilde Z +\lambda) =
(f/\pi^*\theta)(\tilde Z ),\ \hbox{for all}\ \lambda\in\Lambda^\prime,\ \tilde Z\in\C^{g+n}.
\end{equation}
It is easy to compute that 
\[
\theta_{{k_1},\ldots,{k_n}}(z_1,\ldots,z_{g+n})=\exp(2\pi i\sum_{j=1}^n
k_jz_{g+j})\theta(Z+\sum_{j=1}^nk_j\caA(O_{j+1}-O_1))
\]
where $Z=(z_1,\ldots,z_g)$ and the projection $\C^{g+n}\to \C^g$ covering $J^\prime\to J$ has fibre
coordinates $z_{g+1},\ldots,z_{g+n}$. A straightforward but tedious Fourier decomposition argument
shows that all holomorphic functions $f(\tilde Z)$ satisfying (\ref{eq:functional}) must be a linear
combination of the functions above (cf.\ the proof that the classical $\theta$-function spans
$\Gamma(\caT)$ in e.g.\ \cite{GriH}).$\Box$

We want to use this to produce a more explicit expression for the maps
$\psi_{\caL^\prime}:N^\prime\to \Gr$ which appear in theorem \ref{th:Gr}. To this end, fix
$\caL^\prime\in N^\prime$ and set $\caL=\pi(\caL)$. Let us begin by considering, for $L^\prime\in
J_R^\prime$ and $L=\pi(L^\prime)$, the maps
\begin{equation}
\label{eq:CP}
\P\Gamma(\caL\otimes L)\to \CP^n,\quad e_L\mapsto 
[(e_L/(\sigma\otimes s_L))(O_1),\ldots,(e_L/(\sigma\otimes s_L))(O_{n+1})]
\end{equation}
derived from the canonical trivialisation (cf.\ equation (\ref{eq:canon2})). Recall that in this
expression $\sigma$ generates $\Gamma(\caL^\prime)$ and $s_L$ is a trivialising section for
$\caL^\prime\otimes L^\prime$ over $X-\lambda^{-1}(\{0,\infty\})$.
We are only interested in the case where $e_L$ generates $\Gamma(\caL(-\tilde D)\otimes L)$ for some
positive divisor $\tilde D$ of degree $n$. In that case there is a positive divisor $D$ of degree
$g$ such that $\caL(-\tilde D)\simeq \caO_X(D)$. Let $\caT_D$ be the translate of $\caT$ for which
$\caA^*\caT_D\simeq\caO_X(D)$, then the previous proposition holds equally well for $\caT_D$ with
$\theta(Z)$ replaced by the appropriate translate $\theta(Z+\kappa_D)$. 
By (\ref{eq:Theta}) the monomials
$1,x_1,\ldots,x_n$ have images in $\Gamma(\pi^*\caT_D)$ which we will simply call 
$\theta_0(\tilde Z+\kappa_D),\theta_1(\tilde Z+\kappa_D), \ldots ,\theta_n(\tilde Z +\kappa_D)$
(here we lift $\kappa_D\in\C^g$ up to $\C^{g+n}$ using the first $g$ coordinates of the latter).
\begin{prop}
For $e_L\in\Gamma(\caL(-\tilde D)\otimes L)$ as described above the map (\ref{eq:CP}) can also be
written
\begin{equation}
\label{eq:theta}
e_L\mapsto [c_0\theta_0(\tilde Z +\kappa_D),\ldots,c_n\theta_n(\tilde Z +\kappa_D)]
\end{equation}
where $\C^{g+n}\to J^\prime\to J$ maps $\tilde Z\mapsto L$ and 
$c_0,c_1,\ldots,c_n$ are constants which depend only on $\caL^\prime$ and $\tilde D$.
\end{prop}
{\bf Proof.} Fix the base point for $\caA:X\to J$ to be $O_1$. Then $e_L$ can be identified with
$\caA^*\theta(Z+\kappa_D)$ and so $e_L\vert O_k$ is identified with
\[
\theta(\caA(O_k)+Z+\kappa_D)=t_k^{-1*}\theta(Z+\kappa_D),
\]
since $\caA(O_k)=\caA(O_k-O_1)$ with this choice of base point. Now we observe that we have
the tautological sections
\[ 
s_k:J^\prime\to\pi^*J_k, \quad L^\prime\to (L^\prime, s_L\vert O_1\otimes s_L^{-1}\vert O_k).
\]
It follows that
\[
[(e_L\otimes s_L^{-1})\vert O_1,\ldots,(e_L\otimes s_L^{-1})\vert O_{n+1}] =
[\theta_0(\tilde Z +\kappa_D),\ldots,\theta_n(\tilde Z +\kappa_D)].
\]
We obtain (\ref{eq:theta}) by observing that the additional information 
$\sigma\vert O_1,\ldots,\sigma\vert
O_{n+1}$ is independent of $L^\prime$ and only contributes constants $c_0,\ldots,c_n$.
$\Box$
 
Let us apply this to write down the formula for the map $\varphi:\R^{2k}\to\Gr$ with spectral data
$(X,\lambda,\caL)$. For simplicity, set $z=(z_1,\ldots,z_k)$, let $U\in\C^{g+n}$ be the tangent
vector $d\gamma(\partial/\partial z_1,\partial/\partial z_k)$ at $z=0$ and define $U\cdot z=\sum
U_jz_j$. The space $\Gamma(\caL(-D_\infty))$ is spanned by non-zero sections $v_1,\ldots,v_k$ 
for which $v_j$ has divisor of zeroes 
$\hat Q_j=(\lambda)_\infty - Q_j$: these will not be orthogonal. In the linear equivalence 
class for $\caL(-\hat Q_j))$ there is a unique positive divisor $D_j$:
let $\kappa_j\in\C^g$
be such that $\theta(\caA_{P_0}(P)+\kappa_j)$ has divisor of zeroes $D_j$. Then 
\[
\varphi(z) = v_1(z)\wedge \ldots\wedge v_k(z)
\]
where 
\[
v_j(z) = (c_{0j}\theta_0(U\cdot z+\bar U\cdot\bar z +\kappa_j),\ldots,
          c_{nj}\theta_n(U\cdot z+\bar U\cdot\bar z +\kappa_j)).
\]
The constants $c_{ij}$ can be computed using $\theta$-functions but the
formula is not elementary.

\subsection{Algebraic harmonic tori.}

From the results above (see also \cite{McI2}) it is clear that we can obtain
harmonic tori whenever the map $\gamma:\R^{2k}\to J^\prime_R$ has two
independent periods. In general this is a difficult condition to study, since
the image of $\gamma$ need not be algebraic
and in general we must solve some transcendental equations in
the moduli of $(X,\lambda)$. Here we
examine what happens when the image of $\gamma$ is algebraic. For
simplicity we will work with the case of maps into $\CP^n$. We will say a map
$\varphi:\R^2\to\CP^n$ is algebraic if the image of the complexification of
$\gamma$ is an algebraic subgroup of $J^\prime$. We can describe this complex
group homomorphism $\gamma^\C$, as follows. Recall that
to construct maps into $\CP^n$ with spectral data $(X,\lambda)$ we first
single out a zero $P_1$ of $\lambda$ (of at least degree two) 
and the corresponding pole $Q_1=\rho_*P_1$. Let
$\caA^\prime_{P_1},\caA^\prime_{Q_1}:X^\prime\to J^\prime$ be the (rational) Abel maps for these base
points. Then $\gamma^\C:\C^2\to J^\prime$ is the unique homomorphism of complex
groups with the property that
\[
(\partial\gamma^\C/\partial z)(0,0) =
(\partial\caA^\prime_{P_1}/\partial\zeta)(P_1) ;\quad
(\partial\gamma^\C/\partial w)(0,0) =
(\partial\caA^\prime_{Q_1}/\partial\eta)(Q_1),
\]
for coordinates $(z,w)$ on $\C^2$ and local parameters
$\zeta=\sqrt{\lambda}$ about $P_1$ and $\eta=\sqrt{\lambda^{-1}}$ 
about $Q_1$. We may choose these square roots so that $\overline{\rho_*\zeta}=
-\eta^{-1}$. It follows that $\gamma^\C(z,\bar z)$ is the real homomorphism $\gamma(z)$. We will
denote the image of $\gamma^\C$ by $M\subset J^\prime$ and let $M_R$ be the
image of $\gamma$ (which is the connected component of the identity of $M\cap J^\prime_R$).
One knows from \cite{McI2}, p 520, that $M$ is two dimensional unless 
$X$ is the Riemann sphere and $\lambda$ has degree two (for 
that is the only case for which $\dim J^\prime <2$). 
\begin{prop}
If $M\subset J^\prime $ is an algebraic subgroup then it must be one of
three types: (A) $M\simeq \C^\times\times\C^\times$. In this case $X$ is the Riemann sphere and the 
harmonic map $\varphi:M_R\to\P^n$ is $S^1\times S^1$-equivariant. (B) $M$ is a $\C^\times$-extension of an
elliptic curve. In this case $X$ is hyperelliptic and the map $\varphi$ is $S^1$-equivariant. (C) $M$ is
compact i.e.\ an abelian surface. In each case $M_R$ is a real two torus. 
\end{prop}
{\bf Proof.} Set $\hat S=M\cap H_\fo$ and let $S$ be its connected component of the identity. Since $M$ is
algebraic $S$ is a product of $\C^\times$'s and $\hat S/S$ is a finite group. We know immediately from
proposition \ref{pp:equiv} that the map $\varphi:M_R\to\P^n$ must be equivariant for the real group
$\hat S\cap M_R$.
Now set $A=M/S$, then we have exact sequences
\[
1\to S\to M\to A\to 1,\quad 1\to \hat S/S\to A\to \pi(M)\to 1.
\]
Notice that $\pi(M)\subset J$ must be an abelian variety, possibly trivial, and therefore so is 
$A$. Now we have only three
possibilities for the dimension of $S$:

\noindent
(A) $\dim S=2$. It follows that $M=S$ (since $M$ is connected) whence $M\simeq \C^\times\times\C^\times$,
so $\varphi$ is $S^1\times S^1$ equivariant.
But this means $\pi(M)=1$. Now we observe that $M$ is tangent to $\caA^\prime_{P_1+Q_1}
(X^{\prime(2)})\subset J^\prime$ (where
$X^{\prime(2)}$ denotes the symmetric product) and one knows that $d\pi\circ d\caA^\prime=d\caA$. It follows
that $\pi(M)$ is tangent to $\caA_{P_1+Q_1}(X^{(2)})\subset J$. 
Therefore to have $\pi(M)=1$ it must be that $J$ is
the trivial group i.e.\ $X$ is the Riemann sphere. 

\noindent
(B) $\dim S=1$. It follows that $\dim A =1$ i.e.\ it is an elliptic curve so $M$ is a
$\C^\times$-extension and $\varphi$ is $S^1$ equivariant. We also have $\dim\pi(M)=1$. 
But one easily sees, by identifying $\Gamma(\Omega_X)^*$ 
with the tangent space to $J$ at the identity, that the image of $d\caA_{P_1+Q_1}$
is the annihilator of $\Gamma(\Omega_X(-P_1-Q_1)$, so
this must have codimension one in
$\Gamma(\Omega_X)$. It follows that $P_1+Q_1$ is the divisor of zeroes of a rational function on $X$ 
(see e.g.\ \cite{Har} p 341) i.e.\ $X$ is hyperelliptic. 

\noindent
(C) $\dim S=0$. In this case $M$ equals $A$ and is therefore isogenous to $\pi(M)$, 
hence $M$ is an abelian surface.  

In each case $S\cap J^\prime_R$ and $A\cap J^\prime_R$ are real tori, therefore $M_R$ is a real 2-torus.
$\Box$

\noindent
{\em Remarks.} The algebraic tori are a very special class, for in this case the $\theta$-functions
which occur in (\ref{eq:theta}) will reduce to $\theta$-functions for the generalised abelian variety
$M$. In case (A) $M$ is the generalised Jacobian of a rational curve with two nodes, while in case
(B) $M$ is the generalised Jacobian of an elliptic curve with one node (but note that in case
(C) $M$ need not be the Jacobi variety of any curve). However, 
the case of totally equivariant minimal tori in $\CP^n$ (i.e.\ case (A) for conformal maps)
has already been treated by Jensen-Liao \cite{JenL} using a
parameterisation due to Kenmotsu \cite{Ken}. I don't believe that the spectral
curve approach adds much to this. One knows from \cite{McI1} that a harmonic map 
$\varphi:\R^2\to\CP^n$ is totally equivariant if and only if its spectral curve is the 
Riemann sphere and since $J^\prime$ is a product of $\C^\times$'s it follows that $\varphi$ is
algebraic if and only if it is a two torus.
One can write down fairly explicitly the conditions which guarantee that $M$ is algebraic. 
They amount to
a system of linear equations, parameterised by the moduli of the pair $(X,\lambda)$,  which must have
rational solutions (the reader interested in seeing a description of the moduli space $(X,\lambda)$ in
this case is referred to \cite{Tan}). 

The other two cases are more interesting but more enigmatic. Indeed, the only examples I am currently
aware of are the maps $T^2\to\CP^1$ with elliptic spectral curve: these are all Gauss maps of
the Delaunay surfaces of constant mean curvature in Euclidean 3-space. Here $\dim J^\prime = 2$ so
the algebraic condition is automatically satisfied.

\smallskip\noindent
{\small
Dept.\ of Mathematics, University of York, Heslington, York YO10 5DD, U.K.
Email: im7@york.ac.uk
}

\end{document}